\documentclass[reqno]{amsart}

\usepackage[english]{babel}
\usepackage{graphicx}
\usepackage{textcomp}
\usepackage{amsmath}
\newtheorem{thm}{ Theorem}[section]

\begin{document}

\title [Suppression of unbounded gradients ]
{Suppression of unbounded gradients in a SDE associated with the
Burgers equation}

%\title{The non-viscous Burgers equation in a noisy medium: a threshold effect for blow up behaviour}

\author[Albeverio,\,Rozanova]{Sergio Albeverio $^{1}$, Olga Rozanova $^{2}$}

%\author[label1,label2]{}
\address[$^{1}$]{Universit\"{a}t Bonn,
Institut f\"{u}r Angewandte Mathematik, Abteilung f\"{u}r
Stochastik, Wegelerstra\ss e 6, D-53115, Bonn; HCM, SFB611 and IZKS,
Bonn; BiBoS, Bielefeld--Bonn, Germany}
\address[$^{2}$]{Mathematics and Mechanics Faculty, Moscow State University, Moscow
119992, Russia}

\thanks {Supported by Award DFG 436 RUS 113/823/0-1 and the special program of the
Ministry of Education of the Russian Federation "The development of
scientific potential of the Higher School", project 2.1.1/1399}

\email[$^{1}$]{albeverio@uni-bonn.de}
\email[$^{2}$]{rozanova@mech.math.msu.su}



\date{\today}

\begin{abstract}
We consider the Langevin equation describing a non-viscous Burgers
fluid stochastically perturbed by uniform noise. We introduce a
deterministic function that corresponds to the mean of the velocity
when we keep fixed the value of the position. We study
interrelations between this function and the solution of the
non-perturbed Burgers equation. Especially we are interested in the
property of the solution of the latter equation to develop unbounded
gradients within a finite time. We study the question how the
initial distribution of particles for the Langevin equation
influences this blowup phenomenon. We show that for a wide class of
initial data and initial distributions of particles the unbounded
gradients are eliminated. The  case of a linear initial velocity is
particular. We show that if the initial distribution of particles is
uniform, then the mean of the velocity for a given position
coincides with the solution of the Burgers equation and,in
particular, it does not depend on the constant variance of the
stochastic perturbation. Further, for a one space  variable we get
the following result: if the decay rate of the even power-behaved
initial particles distribution at infinity is greater or equal
$|x|^{-2},$ then the blowup is suppressed, otherwise, the blowup
takes place at the same moment of time as in the case of the
non-perturbed Burgers equation.
\end{abstract}

%\pacs{}
\maketitle

\medskip

\section{Introduction}

It is well known that the non-viscous Burgers equation, the simplest
equation that models the nonlinear phenomena in a force free mass
transfer,
$$
u_t+(u,\nabla)\,u=0,\eqno(1.1)
$$
where $u(x,t)=(u_1,...,u_n)$ is a vector-function ${\mathbb
R}^{n+1}\to{\mathbb R}^n,$ before the formation of shocks is
equivalent to the system of ODE
$$
\dot x(t)= u(t,x(t)),\quad \dot  u(t,x(t))=0.\eqno(1.2)
$$
The latter system defines a family of characteristic lines $x=x(t)$
that can be interpreted as the Lagrangian coordinates of the
particles.

Given initial data $$u(x,0)=u_0(x),\eqno(1.3)$$ one can readily get
an implicit solution of (1.1), (1.3), namely, $$u(t,x)=u_0(x-t
u(t,x)).$$ For special classes of initial data we can obtain an
explicit solution. The  simplest case is
$$u_0(x)=\alpha x,\quad \alpha\in{\mathbb R}, \eqno(1.4)$$
where
$$u(t,x)=\frac{\alpha
x}{1+\alpha t} . \eqno(1.5)
$$
Thus, if $\alpha<0,$ the solution develops a singularity at the
origin as $t\to T,\,0<T<\infty,$ where
$$
T=-\frac{1}{\alpha}.\eqno(1.6)$$

 In the present paper we consider a $2\times n$ dimensional
It$\rm\hat{o}$ stochastic differential system of equations,
associated with (1.2), namely
$$
d X_k(t)=U_k(t)\,dt, $$
$$d U_k(t)=\,\sigma \,d(W_k)_t,\quad k=1,..,n,
$$
$$
X(0)=x,\quad U(0)=u,\quad t\ge 0,
$$
where $(X(t),U(t))$ runs the phase space ${\mathbb R}^n\times
{\mathbb R}^n,$ $\sigma>0$  is constant, $(W)_{k,t},\,k=1,...,n,$ is
the $n$ - dimensional Brownian motion.

Our main question is: can a stochastic perturbation  suppress the
appearance of unbounded gradients?

%In the theory of stochastic dynamical systems often consider a
%stochastic perturbation of the velocity, which leads to the
%appearance of a white noise in the second of equations (1.7). The
%problem of solving such equations was investigated in many works
%(see, e.g. \cite{AKl1}, \cite{AK2}, \cite{AK}, \cite{AHZ}, \cite{Risken}). This
%type of stochastic perturbation corresponds to the stochastically
%forced Burgers equation, or in the language of physicists, Burgers
%turbulence. This has been an area of intensive research activity in
%the last decade (see e.g. \cite{Woy}, and for a very recent review
%\cite{Khanin}, and references therein).

The stochastically perturbed Burgers equation and the relative
Langevin equation were treated in many works (e.g.
\cite{Risken},\cite{Sinai}).
 The behavior of the space gradient of the velocity was studied
earlier in other contexts in \cite{Bouchaud}, \cite{Gurarie}, but
this problem is quite different from the problem considered in the
present paper. The analogous problem concerning the behavior of
gradients of solutions to the Burgers equation under other types of
stochastic perturbations was studied in \cite{AR}.

Let us consider the mean of the velocity $U(t)$ at time $t$ when we
keep the value of $X(t)$ at time $t$ fixed but allow $U(t)$ to take
any value it wants, namely
$$
\hat u(t,x)=\frac{\int\limits_{{\mathbb
R}^n}\,u\,P(t,x,u)\,du}{\int\limits_{{\mathbb
R}^n}\,P(t,x,u)\,du},\quad t\ge 0,\,x\in {\mathbb R}^n, \eqno(1.7)
$$
where $ P(t,x,u)$ is the probability density in position and
velocity space, so that ${\int\limits_{{\mathbb R}^n\times {\mathbb
R}^n}}\,P(t,x,u)\,dx\,du=1.$

This function obeys the following Fokker-Planck equation:
$$
\frac{\partial P(t,x,u)}{\partial t
}=\left[-\sum\limits_{k=1}^n\,u_k\,\frac{\partial }{\partial x_k }\,
 + \frac{1}{2}\, \sigma^2
\,\frac{\partial^2 }{\partial u_k^2}\,\right] P(t,x,u), \eqno(1.8)
$$
subject to the initial data
$$P(0,x,u)=P_0(x,u).$$

%Let us denote by $\hat u(t,x)$ the conditional expectation of the
%velocity $U(t)$ at time $t$ given the position $X(t)$ at time $t.$

If we choose
$$
P_0(x,u)=\delta (u-u_0(x))\,f(x)=\prod\limits_{k=1}^{n}\,\delta
(u_k-(u_0(x))_k)\,f(x),\eqno(1.9)$$ with an arbitrary sufficiently
regular $f(x),$ then
$
\hat u(0,x)=u_0(x).
$
The function $f(x)$ has the meaning  of a probability density of the
particle positions in the space at the initial moment of time and
therefore  $f(x)$ has to be chosen such that
$\displaystyle\int\limits_{{\mathbb R}^n}\,f(x)\,dx=1$. If the
latter integral diverges for a certain choice of $f(x),$ we consider
the domain $\Omega_{L}:=[-L,L]^n, \,L>0 $ and the re-normalized
density $f_L(x):=\,\chi(\Omega_{L})\,
f(x)\,\left(\displaystyle\int\limits_{\Omega_{L}}
\,f(x)\,dx\right)^{-1},\,$ where $\chi (\Omega_{L})$ is the
characteristic function of $\Omega_{L}, $ we denote the respective
probability density in velocity and position by $P_L(t,x,u)$ and
modify the definition of $\hat u (t,x)$ as follows:
$$
\hat u(t,x)=\,\lim\limits_{L\to\infty}\,\frac{\int\limits_{{\mathbb
R}^n}\,u\,P_L(t,x,u)\,du}{\int\limits_{{\mathbb
R}^n}\,P_L(t,x,u)\,du},\quad t\ge 0,\,x\in \Omega_L,\eqno(1.10)
$$
provided the limit exists.

%For the sake of simplicity we consider initial data $u_0(x)$ such
%that the determinant of the Jacobi matrix $J(u_0(x))$
%does not
%vanish and

We apply heuristically the Fourier transform in the variables $u$
and $x$ to (1.8), (1.9) to obtain for $\,\tilde P\,= \,\tilde
P(t,\lambda, \xi)\,$
$$
\frac{\partial  \tilde P}{\partial t}\,=\,-\frac{\sigma^2}{2}
\,|\xi|^2 \tilde P\,+\,(\lambda ,\frac{\partial  \tilde P}{\partial
\xi}),\eqno(1.11)$$
%!!!!!!!!!zdes' ispravlena opechatka!!!!
$$
\tilde P(0,\lambda, \xi)\,=\,\int\limits_{{\mathbb R}^n}\,f(s)
e^{-i(\xi, u_0(s))}\,e^{-i(\lambda, s)}\,ds.\eqno(1.12)
$$
Thus, (1.11) and (1.12) give
$$
{\tilde P}(t,\lambda,
\xi)\,=\,e^{-\frac{\sigma^2}{6|\lambda|}(|\xi+\lambda
t|^3-|\xi|^3)}\, \int\limits_{{\mathbb R}^n}\,f(s)
e^{-i((\xi+\lambda t), u_0(s))}\,e^{-i(\lambda, s)}\,ds
\,,\eqno(1.13)$$
$$
{P}(t,x, u)\,=\frac{1}{(2\pi)^{2n}}\,\int\limits_{{\mathbb
R}^{2n}}\,{\tilde P}(t,\lambda, \xi)\,e^{i(\xi, u)}\,e^{i(\lambda,
x)}\,d\lambda\,d \xi\,=
$$
$$=
\left(\frac{\sqrt{3}}{\pi\sigma^2
t^2}\right)^n\,\int\limits_{{\mathbb
R}^n}\,f(s)\,e^{-\frac{2}{\sigma^2\,t^3}\,(3 t^2\,(u, u_0(s))
+t^2\,|u_0(s)-u|^2+3|x-s|^2+3t\,(u+u_0(s),s-x))}\,ds. \eqno(1.14)
$$
Now we substitute (1.14) in (1.8) or (1.10), integrate with respect
to $u$ and get the formula
$$
\hat u(t,x)\,=\,\frac{1}{2t}\,\frac{\int\limits_{{\mathbb
R}^n}\,(-u_0(s)t-3(s-x))\,f(s)\,e^{-\frac{3
|u_0(s)t+(s-x)|^2}{2\sigma^2 t^3}}\,ds}{\int\limits_{{\mathbb
R}^n}\,f(s)\,e^{-\frac{3 |u_0(s)t+(s-x)|^2}{2\sigma^2
t^3}}\,ds},\quad t\ge 0,\,x\in {\mathbb R}^n, \eqno(1.15)
$$
provided all integrals exist.

 Thus, we can compare $\hat u(t,x)$ with
the solution  $u(t,x)$ of the non-viscous Burgers equation (1.1).

%\section{Linear initial velocity}

\section{Exact results}

It is natural to begin with the case where the solution to the
Burgers equation (1.1) can be obtained explicitly. Let us choose
$$u_0(x)=\alpha x,\quad\alpha<0.\eqno(2.1)$$
One can see from (1.5), (1.6) that the gradient of the solution
becomes unbounded as $t\to T.$

If the initial distribution of particles is either uniform or
Gaussian, it is possible to get explicit formulas for $\hat u.$
Namely, for the uniform distribution $f(x)=\rm const$ both integrals
in the numerator and the denominator in (1.15) can be taken and we
get
$$\hat u(t,x)=\frac{\alpha
x}{1+\alpha t},$$ which coincides with (1.5). Therefore,  the
gradient becomes unbounded at $T=-\frac{1}{\alpha}.$ On the
contrary, in the case of a Gaussian distribution,
$f(x)=\left(\frac{r}{\sqrt{\pi}}\right)^n\,e^{-r^2 |x|^2},$  $r>0,$
we get another explicit formula:
$$
\hat u(t,x)=\frac{3(\alpha(\alpha t+1)+r^2 \sigma^2
t^2)}{3(at+1)^2+2 r^2 \sigma^2 t^3}\,x.\eqno(2.2)
$$
One can see that the denominator does not vanish for all positive
$t,$ and at the critical time $T$ we have $\hat
u(t,x)=-\frac{3}{2}\alpha x, $ that is the gradient becomes positive
and tends to zero as $t\to+\infty.$

\section{1D case, specific classes of initial distributions of particles and initial data}

Our main question is how the decay rate of the function $f(x)$ at
infinity relates to the property of $\hat u$ to reproduce the
behavior of the solution of the non-perturbed Burgers equation at
the critical time. For the sake of simplicity we dwell on the case
of an one dimensional space, however the results can be extended to
the higher dimensional space. Let us consider the class of initial
distributions of particles $f(x)$ which are intermediate between
Gaussian and uniform. Our aim is to find a threshold rate of decay
at infinity that still allows to preserve the singularity at the
origin.

We restrict ourselves to the class of smooth distributions $f(x)$
and initial data $u_0(x)$ satisfying the condition
$$
\Big|\int\limits_{{\mathbb R}}\xi^m
\,(u_0(\xi))^l\,f(\xi)\exp\left(-\gamma
\xi^2\right)\,d\xi\Big|<\infty \quad\mbox{for all} \quad
m,l\in{\mathbb N}\cup\{0\},\,\gamma>0.\eqno(3.1)
$$
As a representative of such a class of distributions we can consider
$$f(x)={\rm const}\cdot (1+|x|^2)^{k},\,k\in {\mathbb R}.\eqno(3.2)$$
\begin{thm} Let the initial be $u_0(x)$ be smooth,
 and  for a certain fixed $\beta<0 $ and all $x\in {\mathbb R}$ (except for maybe
  a bounded set)
  $\,|{u_0(x)}-\beta {x}|\ge \gamma>0.$ Moreover, assume that the distribution function $f(x)$ is smooth,
nonnegative, and  the property (3.1) is satisfied. Then the mean
$\hat u (t,x)$ has at the origin $x=0$ at the moment
$t_0=-\frac{1}{\beta},\,\beta<0,$ a bounded derivative
$u'_x(t_0,0).$
%The sign of this derivative is defined by the sign of
%integral in the nominator of (3.3).
\end{thm}
We remark that the initial data with a linear initial profile except
for $\beta=u'_0(x)$ fall into the class of initial data that we have
described above.

\begin{proof}
First of all we perform a change of the time variable. Let
$\epsilon=t+\frac{1}{\beta},\, \beta<0.$ We expand $\hat
u(t(\epsilon),x)$ given by (1.15) into Taylor series at the point
$t=t_0=-\frac{1}{\beta}\,(\epsilon=0),\,x=0,$  taking into account
that condition (3.1) guarantees the convergence of the integrals in
the coefficients of the expansion. This expansion has the form
$$
\hat u(t(\epsilon),x)\sim \frac{1}{2}\,\frac{\int\limits_{\mathbb
R}\left(3 s\beta +u_0(s)\right)\,
f(s)\,e^{\frac{3\beta^3}{2\sigma^2}\left(\frac{
u_0(s)}{\beta}-s\right)^2}\,ds }{\int\limits_{\mathbb R} \,
f(s)\,e^{\frac{3\beta^3}{2\sigma^2}\left(\frac{
u_0(s)}{\beta}-s\right)^2}\,ds }\,+$$
$$
-\frac{3\beta}{2\sigma^2}\left(\frac{\int\limits_{\mathbb
R}\left(\sigma^2-4 \beta^2 s u_0(s)+3s^2\beta^3+\beta
(u_0(s))^2\right)\, f(s)\,e^{\frac{3\beta^3}{2\sigma^2}\left(\frac{
u_0(s)}{\beta}-s\right)^2}\,ds }{\int\limits_{\mathbb R} \,
f(s)\,e^{\frac{3\beta^3}{2\sigma^2}\left(\frac{
u_0(s)}{\beta}-s\right)^2}\,ds }\right. +$$
$$\left.
\frac{\int\limits_{\mathbb R}\left(\beta s - u_0(s)\right)\,
f(s)\,e^{\frac{3\beta^3}{2\sigma^2}\left(\frac{
u_0(s)}{\beta}-s\right)^2}\,ds \,\int\limits_{\mathbb R}\left(3
\beta s - u_0(s)\right)\,
f(s)\,e^{\frac{3\beta^3}{2\sigma^2}\left(\frac{
u_0(s)}{\beta}-s\right)^2}\,ds}{\left(\int\limits_{\mathbb R} \,
f(s)\,e^{\frac{3\beta^3}{2\sigma^2}\left(\frac{
u_0(s)}{\beta}-s\right)^2}\,ds \right)^2} \right)\,x,
%+o(x)O(\epsilon)
\eqno(3.3)
$$
as $x\to 0,\,\epsilon\to 0-$ (where $\,\sim\,$ stands for the
quotient of the left and right sides converging to 1).

The theorem follows immediately from the asymptotics (3.3).
\end{proof}

Let us notice that for even $f(x)$ and odd $u_0(x)$ the expansion
(3.3) is less cumbersome, namely,
$$
\hat u(t(\epsilon),x)\sim
-\frac{3\beta}{2\sigma^2}\frac{\int\limits_0^\infty\left(\sigma^2-4
\beta^2 s u_0(s)+3s^2\beta^3+\beta (u_0(s))^2\right)\,
f(s)\,e^{\frac{3\beta^3}{2\sigma^2}\left(\frac{
u_0(s)}{\beta}-s\right)^2}\,ds }{\int\limits_0^\infty \,
f(s)\,e^{\frac{3\beta^3}{2\sigma^2}\left(\frac{
u_0(s)}{\beta}-s\right)^2}\,ds }\,x,
%+o(x)O(\epsilon)
\eqno(3.3*)
$$
as $x\to 0,\,\epsilon\to 0-.$

It can be readily calculated that if $\beta \to -\infty $ $\,(t_0\to
0),$ then (3.3*) yields  $\hat u(t,x)\sim \alpha x,$ $\,x\to
0,\,\epsilon \to 0-,\,$ where $\alpha=u_x(0)$ (taking account of
$\frac{u_0(\xi)}{\xi}\sim \alpha,\,u'_0(\xi)\sim \alpha,\, \xi\to
0$).

\subsection{Power-behaved distribution}

Let us consider the specific class of even distributions (3.2) and
linear initial data.

The case a linear initial function $u_0(x)=\alpha
x,\,\alpha\ne\beta,$ is particular. Indeed, we have from (3.3*) for
$x\to 0$ and for $t\to t_0=-\frac{1}{\beta},\,\beta<0,$ the
following asymptotic behavior:
$$
\hat u(t,x)\sim \Lambda(\beta)\,x,
$$
with
$$
\Lambda(\beta)=
-\frac{3\beta}{2\sigma^2}\frac{\int\limits_0^\infty\left(\sigma^2+\beta
s^2(\alpha-\beta)(\alpha-3\beta)\right)\,
f(s)\,e^{\frac{3s^2\beta}{2\sigma^2}\left(\beta-\alpha\right)^2}\,ds
}{\int\limits_0^\infty \,
f(s)\,e^{\frac{3s^2\beta}{2\sigma^2}\left(\beta-\alpha\right)^2}\,ds
}.
$$
 We can see that if
$\beta<\alpha$ (before the critical time $T=-\frac{1}{\alpha},$ when
the solution of the non-perturbed Burgers equation blows up) or
$\beta>\alpha$ (after the time $T$) both integrals in (3.3*)
converge and therefore the derivative $\hat u'_x(t,0)$ remains
bounded.
%Moreover, as in the case of bounded initial data the sign of this
%derivative as $\beta \to 0-,$ \,($t_0\to+\infty$) is opposite to the
%sign of $k$ for $k\ne 0$ and it is positive for $k=0.$

Let us consider now the critical moment of time $t=T,$ where
$\beta=\alpha.$ In this case $\,\frac{u_0(x)}{x}= \beta$ identically
and we do not have a multiplier that guarantees the convergence of
integrals of the form
$$
\int\limits_{{\mathbb R}_+}\xi^m f(\xi)\,d\xi\quad\mbox{for all}
\quad m\in{\mathbb N},
$$
which is necessary for the validity of the asymptotics (3.3*).

However, fortunately, due to the relative simplicity of $f(x)$ we
can compute $\hat u(t,x)$ in the vicinity of the origin directly,
using the formula (1.15), which in this case takes the form
$$
\hat u(t,x)\,=\,\frac{1}{2t}\,\frac{\int\limits_{{\mathbb
R}^+}\,(-\alpha s t-3(s-x))\,(1+s^2)^k\,e^{-3 |\alpha s
t+(s-x)|^2}\,ds}{\int\limits_{{\mathbb
R}^+}\,(1+s^2)^k\,e^{-3|\alpha s t+(s-x)|^2}\,ds},\quad t\ge
0,\,x\in {\mathbb R}. \eqno(3.4)
$$
Computations show that for $k\ne \frac {m}{2},\,m\in {\mathbb Z},$
the asymptotic behavior of (3.4) as $x\to 0,\,\epsilon\to 0-,$ where
$\epsilon=t+\frac{1}{\alpha},$ can be expressed through the Gamma
function and the generalized Laguerre functions
$L(\nu_1,\nu_2,\nu_3)$, see \cite{Ryzhik}. It has the form
$$
\hat u (t,x)\sim \frac{F_1(\epsilon,k,\alpha,
\sigma)}{F_2(\epsilon,k,\alpha, \sigma)}\,x,\eqno(3.5)
$$
$$
F_1(\epsilon,k,\alpha, \sigma)=A_1(k)\,\epsilon^{-2k-2}+
o(\epsilon^{-2k-2})+ A_2(k)\,\epsilon^0+o(\epsilon^0),
$$
$$
F_2(\epsilon,k,\alpha,
\sigma)=A_3(k)\,\epsilon^{-2k-1}+o(\epsilon^{-2k-1})+
A_{4}(k)\,\epsilon^0+o(\epsilon^0),
$$
 where the coefficients $A_i(k),\,i=1,..,4, $ are
as follows:
$$
A_1(k)=\frac{\pi^2
2^{k+2}\sigma^{2k}(4k^2-1)}{3^{k}\,|\alpha|^{5k+1}\,\cos{\pi
k}}\,\Gamma(k+1)\,L(k,-k+\frac{1}{2},0),
$$
$$
A_2(k)=\frac{3\sqrt{6}|\pi\alpha|^{\frac{5}{2}}}{2\sigma
(k+1)}\,\tan(\pi k)\,L(\frac{1}{2},k+\frac{1}{2},0),
$$
$$
A_3(k)=\frac{\pi^2 2^{k+1}\sigma^{2k}(2k-1)}{3^k|\alpha|^{5k+1}
(k+1)(k+2)\,\cos{\pi k}}\,\Gamma(k+3)\,L(k,-k+\frac{1}{2},0),
$$
$$
A_4(k)=\frac{\sqrt{6}\pi^2
|\alpha|^{\frac{3}{2}}(2k+3)\Gamma(k+3)}{\sigma
(k+1)(k+2)\Gamma(k+\frac{5}{2})}\,\tan(\pi k).
$$
 Thus, if $k< -1,$ then the leading
term of the numerator and denominator in (3.5) as $\epsilon\to 0-$
is  $A_2\,\epsilon^0$  and (3.5) can be written as
$$
\hat u (t,x)\sim \frac{A_2(k) \epsilon^0 +o(\epsilon^0)}{A_{4}(k)
\epsilon^0 +o(\epsilon^0)} \,x\sim (B_1(k)+o(\epsilon^0))\,x,\,x\to
0,\eqno(3.6)
$$
where $B_1(k)=\frac{A_2(k)}{A_4(k)}.$

 This signifies that the derivative
$\hat u'_x(t,0)$ tends to a finite limit as $\epsilon\to 0-.$

If $-\frac{1}{2}>k>-1,$ then the leading term of the denominator is
$A_4(k)\epsilon^0. $ Otherwise, if $k>-\frac{1}{2},$ then this
leading term is $A_3(k)\epsilon^{-2k-1}.$ Thus we have for
$-\frac{1}{2}>k>-1$
$$
\hat u (t,x)\sim \frac{A_1(k) \epsilon^{-2k-2}
+o(\epsilon^{-2k-2})}{A_4 (k)\epsilon^0 +o(\epsilon^0)}
$$
and  $$\hat u'_x(t,0)\sim B_2
(k)\cdot\frac{1}{\epsilon^{2k+2}},\quad B_2(k)=
 \frac{A_1}{A_4},\qquad
x\to 0,\,\epsilon\to 0-.\eqno(3.7)$$ At last for $k>-\frac{1}{2}$ we
have
$$
\hat u (t,x)\sim \frac{A_1(k) \epsilon^{-2k-2}
+o(\epsilon^{-2k-2})}{A_3(k) \epsilon^{-2k-1} +o(\epsilon^{-2k-1})}
\,x,\quad x\to 0,\,\epsilon\to 0-,
$$
and $$\hat u'_x(t,0)\sim B_3(k) \cdot \epsilon^{-1},\quad
B_3(k)=\frac{A_1(k)}{A_3(k)}={2k+1}.\eqno(3.8)$$ If $k\in {\mathbb
Z},$ then the numerator and the denominator in the leading term in
the expansion of (3.5) as $x\to 0$ are expressed either through
rational functions ($k\ge 0$) or through a Gaussian distribution
function  ($k< 0$). For $k= \frac {2l+1}{2},\,l\in {\mathbb Z},$ the
coefficient of the leading term is expressed through a fraction of
series consisting of the digamma functions. Anyway, the asymptotics
(3.6) takes place also for $k= \frac {l}{2},\,l\in {\mathbb
Z},\,k\ne-\frac{1}{2},$ and it can be found also as a limit
$\kappa\to k.$ For $k<-1$ the function $\hat u(t,x)$ behaves as in
(3.6), where the coefficient $B_1(k)$ can be calculated either
independently or as  $\lim\limits_{\kappa\to
k}\frac{A_2(\kappa)}{A_4(\kappa)}.$ Since for $k=-1$ the degrees in
$\epsilon^{-2k-2}$ and $\epsilon^0$ coincide, then
$$
\hat u (t,x)\sim \lim\limits_{\kappa\to
-1}\frac{(A_1(\kappa)+A_2(\kappa)) \epsilon^0
+o(\epsilon^0)}{A_{4}(\kappa) \epsilon^0 +o(\epsilon^0)} \,x\sim
(B_4+o(\epsilon^0))\,x,\,x\to 0,\,\epsilon\to 0-,
$$
where
$B_4=\lim\limits_{\kappa\to-1}\left(B_1(\kappa)+\frac{A_1(k)}{A_4(k)}\right)=\frac{3|\alpha|}{2}-\frac{\sqrt{6}|\alpha|^{\frac{5}{2}}}{\sigma\sqrt{\pi}}.$
For $k\ge 0$ the function $\hat u (t,x)$ has the asymptotics (3.8)
with the same value $B_3(k).$ An exceptional case is
$k=-\frac{1}{2},$ where
$$F_1(\epsilon,-1/2,\alpha, \sigma)=\bar A_1\,{\epsilon^{-1}}+
o\left(\epsilon^{-1}\right), \quad\bar
A_1=\lim\limits_{k\to-1/2}\,A_1=\frac{4\sqrt{6}|\pi\alpha|^{\frac{3}{2}}}{\sigma},$$
$$F_2(\epsilon,-1/2,\alpha, \sigma)=A_5\,\ln(-\epsilon)+
o\left({\ln(-\epsilon)}\right), \quad A_5=-\bar A_1,\,\epsilon\to
0-.$$ Thus, for $k=-\frac{1}{2}$ we have
$$
\hat u_x(t,0)\sim -\,\frac{1}{\epsilon \,\ln(-\epsilon)} +
o\left(\frac{1}{\epsilon \,\ln(-\epsilon)}\right),\,\epsilon\to
0-.\eqno(3.9)
$$
 The following theorem
summarizes our results:
\begin{thm}
Assume that in the case  $n=1$  the initial distribution function is
$f(x)={\rm const}\cdot (1+|x|^2)^{k},\,k\in {\mathbb R},$ and the
initial velocity has  the form  $u_0(x)=\alpha x,\, \alpha<0.$ Then
the space derivative of the mean $\hat u(t,x)$ at the origin $x=0$
is bounded for all $t>0$ except for the critical time
$T=-\frac{1}{\alpha}.$ At the critical time the behavior of the
derivative depends on $k.$ Namely,
 for $k>-1$
the mean $\hat u(t,x)$ keeps the property of solutions to the
non-perturbed Burgers equation  to blow up at the critical time $T$
at $x=0.$ The rates of the blowup for $-\frac{1}{2}>k>-1,$
$k>-\frac{1}{2}$  and $k=-\frac{1}{2}$ are indicated in (3.7), (3.8)
and (3.9), respectively. Otherwise, if $k\le-1,$ the derivative
$\hat u'_x(t,0)$ at the critical time remains bounded, i.e the
singularity disappears.
\end{thm}

%\vskip1cm

\section{Pressureless gas dynamics model and a limit case at vanishing noise}

Let us consider the pair $(\rho_t, u(\cdot,t)),t\ge 0,$ where
$\rho_t$ is the probability distribution of the random variable
$X_t,$ governed by the SDE $X_t=X_0+\int_0^t\,{\mathbb
E}[u_0(X_0)|X_s]\,ds,\,t\ge 0,$ with a given random variable $X_0$
with values in $\mathbb R$ and function $u_0(x),\,$ and
$u(t,x)={\mathbb E}[u_0(X_0)|X_t=x].$ According to \cite{Dermoune},
$(\rho_t, u(\cdot,t),t\ge 0)$ is a weak solution to the pressureless
gas dynamics model
$$\partial_t \rho + \partial_x (\rho u)=0,\qquad  \partial_t (\rho
u)+\partial_x(\rho u^2)=0.\eqno(5.1)$$ Therefore it is natural to
expect that the limit as $\sigma\to 0$ of the mean $\hat u$ (denoted
by $v(t,x)$) takes part in the solution to (5.1).

For smooth $u_0(x)$ and $f(x)$ this can be readily shown. First of
all, let us introduce the function $\rho(t,x)=\int\limits_{{\mathbb
R}^n}\, P(t,x,u)\,du $ and notice that it satisfies the continuity
equation $\,\partial_t \rho + \partial_x (\rho \hat u)=0. $ Further,
we check that  the function $v(t,x)$ satisfies the Burgers equation
(1.1). Indeed, $\frac{t^{3/2}}{\sigma\sqrt{6\pi}}\,\exp\left(
-\frac{3(u_0(s)t+(s-x))^2}{2\sigma^2 t^3}\right)\to \delta
(s-s(t,x)),$ as $\sigma\to 0$ in ${\mathcal D}',$ where $s(t,x)$ is
a solution to equation $u_0(s)+\frac{s-x}{t}=0,$ given in implicit
form. This function exists and it is differentiable provided $t \ne
-\frac{1}{u'_0(x)}$ (at this moment of time the solution to the
Burgers equation blows up).  Thus, $\hat u(t,x)\to u_0(s(t,x))$ as
$\sigma\to 0.$ Now it is sufficient to substitute $u_0(s(t,x))$ into
(1.1) and compute the derivatives of $s(t,x)$ by means of the
implicit function theorem.

 Thus, for smooth initial data $(f(x), u_0(x))$ the pair $(\rho, \hat u)$ is a solution to the system
 $$\partial_t \rho + \partial_x (\rho \hat u)=0,\qquad  \partial_t (\rho
\hat u)+\partial_x(\rho \hat u^2)= \Lambda,\qquad
\Lambda=-\int\limits_{\mathbb R} \,P_x(t,x,u)\,(u-\hat u)^2\,du,$$
where  $\,\Lambda\to 0$ as $\sigma\to 0.$ The integral relaxation
term $\Lambda$ can be used instead of the traditional viscosity
\cite{TZZ}.
%\medskip

It is interesting to consider the Fokker-Plank equation (1.8) as a
kinetic equation
$$
\frac{\partial P(t,x,u)}{\partial t
}+\sum\limits_{k=1}^n\,\left(u_k\,\frac{\partial P(t,x,u)}{\partial
x_k }\,
 +
\,\frac{\partial \,\dot u_k  P(t,x,u)}{\partial u_k}\right)\,=0
$$
(e.g.\cite{Struchtrup}), where the acceleration $\dot u$ of the
particles is due to external forces and the interaction forces with
other particles. It can be readily calculated that in our case $\dot
u=\frac{2}{t}(\hat u-u).$

The authors would like to thank an anonymous referee for helpful
suggestions.

%\newpage


\begin{thebibliography}{99}
\bibitem{Risken}H. Risken, \, {\it The Fokker-Planck
Equation Methods of Solution and Applications,} Second Edition,
Springer-Verlag, 1989.

\bibitem{Sinai} Weinan E, K. M. Khanin, A. E. Mazel, Ya. G. Sinai,\,{\it Invariant measures
for Burgers equation with stochastic forcing,}    Ann. of Math. (2)
151, no. 3, 877-960 (2000).
\bibitem{Bouchaud}J.P.Bouchaud, M.M\'ezard, {\it Velocity fluctuations in forced
Burgers turbulence,} Phys.Rev. E{\bf 54},\, 5116 (1996).
\bibitem{Gurarie}V.Gurarie, \, {\it Burgers equations revisited,}
arXiv:nlin/0307033v1\, [nln.CD]\, (2003).
\bibitem{AR} S.Albeverio, O.Rozanova, \, {\it The non-viscous Burgers equation associated with random positions
in coordinate space: a threshold  for blow up behaviour}, to appear
in $\rm M^3AS$\, (Mathematical methods and models in applied
science)(2009), arXiv:0708.2320v2 [math.AP].
\bibitem{Ryzhik} I.S. Gradshteyn, I.M. Ryzhik,
Table of integrals, series, and products.  6th ed.  San Diego, CA:
Academic Press, 2000.
\bibitem{Dermoune} A.Dermoune, \,{\it Probabilistic interpretation for system of conservation law
arising in advection particle dynamics} C.R.Acad.Sci.Paris, {\bf
326}, Serie I, 595-599 (1998)
%\bibitem{BG} Y.Br\'enier, E.Grenier. {\it Sticky particles and scalar conservation
%laws,} SIAM J.Num.Anal. {\bf 35}, 2317-2328 (1998)
\bibitem{TZZ} D.Tan, T.Zhang, Y.Zheng,\,{\it Delta-shock waves as limits of vanishing viscocity for
hyperbolic system of conservation laws,} J.Diff.Equat., {\bf 112}
(1994), 1-32.
\bibitem{Struchtrup}H.Struchtrup, {\it Macroscopic
Transport Equations for Rarefed Gas Flows: Approximation Methods in
Kinetic Theory,} Springer-Verlag Berlin Heidelberg, 2005.
\end{thebibliography}
\end{document}